\documentclass[11pt,a4paper]{article}
\usepackage[utf8]{inputenc}
\usepackage[T1]{fontenc}
\usepackage[english]{babel}
\usepackage{amsmath,amssymb, amsthm}
\usepackage{mathtools}
\usepackage{lmodern}
\usepackage[margin=2.5cm]{geometry}
\usepackage[colorlinks=true,linkcolor=blue,urlcolor=blue,citecolor=blue]{hyperref}
\usepackage{graphicx}
\usepackage{verbatim}
\usepackage{booktabs}
\usepackage{enumitem}
\newcommand{\oui}{\textbf{Yes}}
\newcommand{\non}{\textbf{Not in general}}
\newcommand{\R}{\mathbb{R}}
\newcommand{\E}{\mathbb{E}}
\newcommand{\C}{\mathbb{C}}

\DeclareMathOperator{\tr}{tr}
\DeclareMathOperator{\diag}{diag}
\DeclareMathOperator{\Vect}{span}
\DeclareMathOperator{\End}{End}

\DeclareMathOperator{\rank}{rank}
\DeclareMathOperator{\Sym}{Sym}

\renewcommand{\Im}{\operatorname{Im}}

\theoremstyle{definition}
\newtheorem{theorem}{Theorem}
\newtheorem{proposition}{Proposition}

\newtheorem{remark}{Remark}
\newtheorem{definition}{Definition}
\newtheorem{lemma}{Lemma}
\newtheorem{corollary}{Corollary}
\title{On a conjecture of \'{E}lie Cartan concerning the Gauss equation\\
 of a flat submanifold}
\author{{Jean-Marie Morvan}\\
Lyon (France)
\\\ \texttt{morvan1704@gmail.com}}
\date{}

\begin{document}
\maketitle
\begin{abstract}
We study exteriorly orthogonal (or flat) symmetric bilinear forms on $\R^n$ with values in $\R^p$, continuing the work initiated by \'{E}lie Cartan in 1919 and 1920. In particular, we describe the second fundamental form of an $n$-dimensional flat manifold isometrically immersed in $\E^{n+p}$. Depending on the values of $n$ and $p$, Cartan's results either can be extended to this setting or fail.
\end{abstract}
\section{Introduction}
This work is motivated by the study of the Gauss equation for an isometric immersion of a Riemannian manifold $M$ into a Euclidean space. The curvature tensor $R$ of $M$ and its second fundamental form $h$ are related by the equation
\begin{equation} \label{E105}
\langle R(X,Y)Z,W\rangle=\langle h(X,Z),h(Y,W)\rangle-\langle h(X,W),h(Y,Z)\rangle.
\end{equation}
Solving this equation amounts to determining all possible $h$ for a given $R$. In other words, if $M$ is known to be isometrically immersed in a Euclidean space, what can be deduced about its second fundamental form, and hence about its ``shape''? The simplest situation is the case where $M$ is flat ($R=0$). Equation~\ref{E105}, usually called the \emph{Gauss equation}, then becomes
\begin{equation} \label{E11}
\langle h(X,Z),h(Y,W)\rangle-\langle h(X,W),h(Y,Z)\rangle=0.
\end{equation}
A partial solution of this linear equation on the tangent bundle of the submanifold is part of \'{E}lie Cartan's foundational work in two papers published in 1919 and 1920 in the \emph{Bulletin de la Soci\'et\'e Math\'ematique de France} \cite{Cartan1919,Cartan1920}. It is a first step towards classifying certain immersions of manifolds of constant curvature into Euclidean spaces or other spaces of constant curvature. Cartan gives the following definition (p.~152; translated here):

\vspace{10pt}

{\it Given $r$ quadratic forms $\Phi_1,\ldots,\Phi_r$ in $n$ variables $x_1,\ldots,x_n$, these forms are said to be} exteriorly orthogonal {\it if the following $\frac{n(n-1)}2$ relations hold:}
$$[\frac{\partial \Phi_1}{\partial x_i}, \frac{\partial \Phi_1}{\partial x_j}] 
+ [\frac{\partial \Phi_2}{\partial x_i}, \frac{\partial \Phi_2}{\partial x_j}]
+ \ldots + [\frac{\partial \Phi_r}{\partial x_i}, \frac{\partial \Phi_r}{\partial x_j}]=0,
 \forall i,j\in\{1,2,\ldots,n\}.$$

\vspace{10pt} 

This definition is followed by two results. The first is the following theorem:

\vspace{10pt}
{\it Given $r$ quadratic forms $\Phi_1,\ldots,\Phi_r$ in $n$ variables that cannot be expressed in terms of fewer than $r$ variables, they are linearly independent and can be obtained, by an orthogonal substitution, from the squares of $r$ linearly independent linear forms.}

\vspace{10pt}

The second is restricted to low dimensions:
\vspace{10pt}

{\it If $n\leq 3$, an orthogonal substitution can always be made so that the forms $\Phi_1,\ldots,\Phi_r$ are perfect squares.}
\vspace{10pt}

Finally, at the end of the proof, Cartan adds the following sentence without further explanation, which we regard as a conjecture:

\begin{center}
{\it The theorem no longer holds for $n>3$.}
\end{center}

Our aim is to understand this last assertion, first by systematically studying the case $n=4$ for the various values of $r$, and then by considering arbitrary $n$ with $r=n+1$. The results vary substantially with $n$ and $r$. The case $n=4$, $r=5$ is similar to Cartan's second theorem, whereas the cases $n=4$, $r=6$ or $7$ are entirely different from those he considered. We have not succeeded in resolving the case $n=4$, $r=8$. The cases $n=4$, $r=9$ or $10$ yield results slightly weaker than Theorems~\ref{T1} and~\ref{T2}. Finally, for arbitrary $n>4$ and $r=n+1$, an additional hypothesis is needed to obtain an analogue of Theorem~\ref{T2}.

\vspace{10pt}

\noindent\textbf{Outline of the paper.}
Section~\ref{P1111} reformulates the definitions in modern terms and records elementary properties of exteriorly orthogonal (or flat) forms. Section~\ref{P1112} states the two results of Cartan mentioned above and outlines their proofs, following both \cite{Cartan1919} and subsequent treatments that expand or simplify the arguments (see, for instance, \cite{Moore1972,Moore1977Positive,DajczerTojeiro2019}). Section~\ref{Par22} introduces our main tool for classifying certain flat forms in higher dimensions. Section~\ref{Par5} systematically studies Cartan's conjecture for $n=4$, considering each value of $r$ from $5$ to $10$. Section~\ref{Par606} gives a solution in the case $r=n+1$. After some perspectives in Section~\ref{Persp}, we add an appendix about semidefinite endomorphisms of $\Sym(\R^n)$, involved in the proofs of our main results.

\section{Flat forms}\label{P1111}
Throughout, $\Sym(V)$ denotes the vector space of symmetric bilinear forms on an $n$-dimensional real vector space $V$, where $n\geq1$; its dimension is $\frac{n(n+1)}2$. Whenever convenient, we identify a symmetric bilinear form with its associated quadratic form and with the corresponding linear map from $V$ to its dual $V^*$. We generally identify $V$ with $\R^n$ after choosing a basis, and $\Sym(\R^n)$ with the space of real symmetric $n\times n$ matrices. When useful, we equip $\R^n$ with its standard inner product $\langle\cdot,\cdot\rangle$, thereby obtaining the Euclidean space $\E^n$. We equip
$\Sym(\R^n)$ with the
Hilbert--Schmidt inner product
$\langle M,N\rangle=\tr(MN)
=\sum_{i,j=1}^n M_{ij}N_{ij}$,
for all $M,N\in\Sym(\R^n)$.

\subsection{Definitions and elementary properties}\label{Par21}
\begin{definition}\label{D1}
Let $V$ be a real vector space of finite dimension $n\geq1$, and let $\{A^1,\ldots,A^p\}\subset\Sym(V)$ be a family of real-valued symmetric bilinear forms on $V$. This family is said to be \emph{exteriorly orthogonal}, or \emph{flat}, if
\begin{equation} \label{E3}
\sum_{{\alpha}=1}^p
\left[
A^{\alpha}(x,z)A^{\alpha}(y, w)
-A^{\alpha}(x,w)A^{\alpha}(y,z)
\right]=0,
\end{equation}
for all $x,y,z,w\in V$.
\end{definition}

Equation~\ref{E3} can be reformulated as follows:
consider the vector-valued symmetric bilinear form ${\cal A}:V\times V\to\E^p$ whose components in an orthonormal basis of $\E^p$ are $(A^1,\ldots,A^p)$. Equation~\ref{E3} then reads
\begin{equation}\label{E4}
\langle{\cal A}(x,z), {\cal A}(y,w)\rangle -\langle{\cal A}(x,w), {\cal A}(y,z)\rangle=0,
\end{equation}
for all $x,y,z,w\in V$. In this case, we also say that $\cal A$ is \emph{exteriorly orthogonal}, or \emph{flat}.
\vspace{10pt}
Let us clarify the notations. If ${\cal A}=(A^1,\ldots,A^p): V \times V \to \E^p$ is any symmetric bilinear map (without any assumption of flatness), we can define three different objects~: the vector subspace $\Vect(\Im {\cal A})$ of $\R^p$ spanned by the image $\Im {\cal A}$ of $\cal A$, the vector subspace  $\Vect (A^1,\ldots,A^p)$ of $\Sym(V)$ spanned by  $A^1,\ldots,A^p$, and the linearisation $L_{\mathcal A}$ defined on the symmetric square of $V$ with values in $\E^p$.  One has 
\begin{equation} \label{E752}
\dim \Vect(\Im {\cal A}) = \dim \Vect(A^1,\ldots,A^p) = \rank L_{\cal A}.
\end{equation}
So we will call the integer defined by equation \ref{E752} the {\it rank of} ${\cal A}$. We will denote by $F$ the subspace  $\Vect (A^1,\ldots,A^p)$ and by $r$ the rank of ${\cal A}$.

 The following observations, which appear in \cite{Cartan1919}, will be used repeatedly.
\begin{enumerate}
\item \label{Par2111}
Obviously,  $r\leq\frac{n(n+1)}2$, independently of whether $\cal A$ is flat.
\item \label{Par2112}
We can always choose a basis $(\epsilon_1,\ldots,\epsilon_p)$ of $\E^p$ such that
$${\cal A}(x,y)= \sum_{\alpha = 1}^{r} A^\alpha (x,y)\epsilon^{\alpha}, \forall x,y \in V,$$
thereby restricting the codomain of $\cal A$ to $\Im {\cal A}$ without loss of generality.

\vspace{10pt}

\noindent Now, for simplicity, we put $V = \R^n$, and suppose that ${\cal A}=(A^1,\ldots,A^p):\R^n\times\R^n\to\E^p$ is flat.
\item \label{Par211}
In any basis $(e_1,\ldots,e_n)$ of $\R^n$ and any orthonormal basis $(\nu_1,\ldots,\nu_p)$ of $\E^p$, writing
${\cal A}^{\alpha}_{ij}=\langle {\cal A}(e_i,e_j),\nu_{\alpha} \rangle,$
the equation becomes
\begin{equation} \label{E5}
\sum_{{\alpha}=1}^p
\left(A^{\alpha}_{ik}A^{\alpha}_{jl}
-A^{\alpha}_{il}A^{\alpha}_{jk}\right)=0.
\end{equation}
We may also choose $(\nu_1,\ldots,\nu_p)$ so that its first $r$ vectors span $\Vect(\Im\mathcal A)$. Equation~\ref{E5} then becomes
\begin{equation} \label{E61}
\sum_{{\alpha}=1}^r
\left(A^{\alpha}_{ik}A^{\alpha}_{jl}
-A^{\alpha}_{il}A^{\alpha}_{jk}\right)=0,
\end{equation}
by choosing an adapted orthonormal basis and using observation~\ref{Par2111}.
We shall refer to any equation of the form~\ref{E3}, \ref{E4}, \ref{E5}, or~\ref{E61} as the \emph{Gauss equation for zero curvature}, or simply the \emph{Gauss equation} when no confusion is possible.
\item \label{Par212}
If $p=1$, $\cal A$ consists of a single form
$A^1$. Equation~\ref{E4} shows that $A^1$ is flat if and only if its rank is at most $1$.
\item\label{Par213}
It follows immediately that if ${\cal A}=(A^1,\ldots,A^p)$ is flat and $q$ of the forms $A^\alpha$ (say, the last $q$) have rank at most $1$, then equation~\ref{E5} reduces to
\begin{equation}\label{E6}
\sum_{{\alpha}=1}^{p-q} \left(A^{\alpha}_{ik}A^{\alpha}_{jl}
-A^{\alpha}_{il}A^{\alpha}_{jk}\right)=0.
\end{equation}
\item \label{Par214}
Equation~\ref{E4} implies that the quadrilinear form $T_{\cal A}:\R^n\times\cdots\times\R^n\to\R$ defined by
\begin{equation}
T_{\cal A}(x,y,z,t)= \langle{\cal A}(x,y),{\cal A}(z,t)\rangle
\end{equation}
is totally symmetric.
In other words, the tensor $T_{\cal A}=\sum_{\alpha=1}^p A^\alpha\otimes A^\alpha$, whose components are $T_{{\cal A}_{ijkl}}=\sum_{\alpha=1}^p A^\alpha_{ij}A^\alpha_{kl}$, is totally symmetric in $i,j,k,l$:
\begin{equation} \label{E9}
T_{{\cal A}_{i_1i_2i_3i_4}}
=T_{{\cal A}_{i_{\sigma(1)}i_{\sigma(2)}i_{\sigma(3)}i_{\sigma(4)}}}, 
\end{equation}
for every permutation $\sigma$.
(Conversely, if $T_{\cal A}$ is totally symmetric, equation~\ref{E4} holds.)
\item \label {Par215}
In the course of proving his first theorem, Cartan establishes the following result (see also \cite{ChernKuiper1952}):

\begin{proposition} \label{P1}
If ${\cal A}=(A^1, \ldots, A^{p}) : \R^n \times \R^n \to \R^p$ is a symmetric bilinear flat form, then 
\begin{equation} \label{E7}
\dim \ker {\cal A} \geq n -\dim \Vect(A^1, \ldots,A^{p}).
\end{equation}
\end{proposition}

\item \label{Par216}
Consider the symmetric bilinear form ${\bf C}:\R^n\times\R^n\to\R$ associated with $T$, obtained by contracting its first two (or last two) indices:

\begin{equation}
{\bf C}_{ij}
=\sum_{\alpha=1}^p\sum_{k=1}^n A^\alpha_{ik}A^\alpha_{kj}
=\sum_{\alpha} \sum_k A_{ij}^{\alpha} A_{kk}^{\alpha}=\sum_{\alpha=1}^p(\operatorname{tr}A^\alpha)A^\alpha_{ij}.
\end{equation}
Then,
${\bf C} \in \Vect(A^1,\ldots ,A^p)$. Moreover,
\begin{equation}
{\bf C}_{ii}=\sum_{\alpha} \sum_k A_{ik}^{\alpha} A_{ki}^{\alpha} \geq 0.
\end{equation}
Suppose now that the codomain is $\E^p$ endowed with its standard scalar product, we get
\begin{equation} \label{E10}
{\bf C}(x,x)=\tr \langle ({\cal A}(.,.),{\cal A}(x,x) \rangle =\langle H,{\cal A}(x,x) \rangle= ||{\cal A}(x,.)||^2 \geq 0, \forall x \in \R^n,
\end{equation}
where $H$ is the trace of $\cal A$, and the norm is the Hilbert--Schmidt norm.
This yields the following proposition.
\begin{proposition}\label{P2}
Let ${\cal A}=(A^1, \ldots, A^{p}) : \R^n \times \R^n \to \E^p$ be a symmetric bilinear flat form with codomain $\E^p$ endowed with  an orthonormal basis. Then, the form ${\bf C}$ is positive semi-definite and $\ker {\bf C}= \ker {\cal A}$. Moreover, ${\bf C}$ is  positive definite if and only if $\ker\mathcal A=\{0\}$.
\end{proposition}
\begin{remark} In submanifold geometry, equation~\ref{E10} means that the $1$-form dual to the mean curvature vector $H$, under the duality induced by the inner product on $\E^p$, is positive semidefinite, and positive definite if the kernel of the second fundamental form is $\{0\}$.
\end{remark}
\end{enumerate}
\section{Two fundamental results of \'{E}. Cartan}\label{P1112}
\subsection{A classification when $\ker\mathcal A=\{0\}$}
A fundamental theorem of Cartan \cite{Cartan1919} classifies flat symmetric bilinear forms ${\cal A}:\R^n\times\R^n\to\R^p$ with $p\leq n$. It is stated for $p=n$, but Proposition~\ref{P1} immediately extends it to $p\leq n$. Indeed, if $\ker\mathcal A\neq\{0\}$, one simply passes to the quotient $\R^n/\ker\mathcal A$ and considers the induced form $\widetilde{\mathcal A}$. Flatness is preserved, yielding analogous results on $\R^{n-\dim\ker\mathcal A}$. This theorem can be stated as follows~:
\begin{theorem} \label{T1}
Let ${\cal A}=(A^1,\ldots,A^n):\R^n\times\R^n\to\R^n$ be a flat symmetric bilinear form with $\ker\mathcal A=\{0\}$. Then, there exist an orthogonal $n\times n$ matrix $(a^\alpha_\beta)$ and $n$ linear forms $\omega^1,\ldots,\omega^n$ such that, for every $\alpha\in\{1,\ldots,n\}$,
$$A^{\alpha} = \sum_{\beta=1}^{n}a_{\beta}^{{\alpha}}\omega^{\beta} \otimes \omega^{\beta}.$$
In particular, $\cal A$ is diagonalizable.
Moreover, the forms $\omega^\beta$ are uniquely determined up to change of sign and  permutation.
\end{theorem}
\noindent The following corollary is immediate.
\begin{corollary}\label{C1}
Let ${\cal A}:\R^n\times\R^n\to\E^n$ be flat, with $\ker\mathcal A=\{0\}$. Then there exist an orthonormal basis $(\epsilon_1,\ldots,\epsilon_n)$ of $\E^n$ and a basis $(\omega^1,\ldots,\omega^n)$ of $(\E^n)^*$ such that
\begin{equation} \label{138 } {\cal A} = \sum_{{\alpha}=1}^{n}\omega^{{\alpha}}\otimes \omega^{{\alpha}}\otimes \epsilon_{{\alpha}}.
\end{equation}
\end{corollary}
The original proof of Theorem~\ref{T1} appears in \cite{Cartan1919}. It has subsequently been revisited, expanded, and extended by several authors \cite{Moore1972,Moore1977Positive,DajczerTojeiro2019,Chen1973}, to whom we refer for details.

\vspace{10pt}

\noindent\textbf{Outline of the proof of Theorem~\ref{T1}.} Here are the main ideas of a simple proof following \cite{Moore1972}. Identify each $A^\alpha$ with a linear map from $\R^n$ to $(\R^n)^*$, and choose $x\in\R^n$ so that the dimension $p$ of the span of $\{A^1(x),\ldots,A^n(x)\}$ is maximal. Suppose that $p<n$. After changing the basis of $(\R^n)^*$, we may assume that $A^{p+1}(x),\ldots,A^n(x)$ vanish. By~\ref{Par213}, the Gauss equation becomes
$$\sum_{{\alpha} = 1}^{p}A^{\alpha}(x) \wedge A^{\alpha}(y)=0,$$
for every $y\in\R^n$. Cartan's lemma then implies that each $A^\alpha(y)$ is a linear combination of the $A^\beta(x)$, $1\leq\beta\leq p$. Suppose that $p<n$. There exist $\alpha\geq p+1$ and $z\neq0$ such that $A^\alpha(x,z)=0$. But there also exist $\alpha\in\{p+1,\ldots,n\}$ and $v\in\R^n$ such that $A^\alpha(z,v)\neq0$. We can then find real numbers $s,t$ such that the span of $A^1(sx+tv),\ldots,A^n(sx+tv)$ has dimension strictly greater than $p$, a contradiction. Thus $p=n$.
To determine the orthogonal matrix $(a^\alpha_\beta)$, choose a basis $(e_1,\ldots,e_n)$ of $\R^n$ with $e_1=x$, and set $\theta^\lambda=A^\lambda(e_1,\cdot)$. Applying Cartan's lemma once more gives
$$A^\lambda (e_i,.)= \sum_\mu C(i)^\lambda_\mu \theta^\mu.$$
The Gauss equation implies that the $C(i)$ commute. Hence they are simultaneously diagonalizable. Let $(a^\alpha_\beta)$ be an orthogonal matrix that diagonalizes them. The forms $\psi^\alpha=\sum_\beta a^\alpha_\beta A^\beta$ have rank $1$, and there exist linear forms $\omega^\alpha$ such that $\psi^\alpha=\omega^\alpha\otimes\omega^\alpha$. We finally obtain
$$A^\alpha = \sum_\beta a_\beta^\alpha  \omega^\beta \otimes \omega^\beta.$$
\noindent\textbf{Outline of the proof of Corollary~\ref{C1}.}
Choose an orthonormal basis $(\eta_1,\ldots,\eta_n)$ of $\E^n$ and set, for all $x,y\in\R^n$,
$
\Phi^\lambda(x,y)=\langle {\cal A}(x,y),\epsilon_\lambda\rangle.
$
Thus,
$$
{\cal A}(x,y)=\sum_{\lambda=1}^n\Phi^\lambda(x,y)\eta_\lambda.
$$
The flatness of $\cal A$ reads
$$\sum_{\lambda=1}^n
\bigl[\Phi^\lambda(x,y)\Phi^\lambda(z,t)
-\Phi^\lambda(x,t)\Phi^\lambda(z,y)\bigr]=0.
$$
The hypotheses of Theorem~\ref{T1} are therefore
satisfied. It yields an orthogonal matrix
$(a^\alpha_\lambda)$ and linear forms $\omega^\alpha$ such that
$$\sum_{\lambda=1}^n a^{\alpha}_{\lambda}\Phi^\lambda
=\omega^{\alpha}\otimes\omega^{\alpha}.
$$
Define
$
\epsilon_{\alpha}=\sum_{\lambda=1}^na_{{\alpha}\lambda}\eta_\lambda.
$
The orthogonality of $(a_\alpha^\lambda)$ implies
$
\langle\epsilon_{\alpha},\epsilon_\beta\rangle
=\sum_{\lambda=1}^n a_{\alpha}^{\lambda}a_{\beta}^{\lambda}
=\delta_{\alpha \beta}.
$
Moreover,
$
\langle{\cal A}(x,y),\epsilon_{\alpha}\rangle
=\omega^{\alpha}(x)\omega^{\alpha}(y).
$
Consequently,
\begin{equation}  {\cal A} = \sum_{\alpha=1}^{n}\omega^{\alpha}\otimes \omega^{\alpha}\otimes \epsilon_{\alpha}.
\end{equation}
Expanding $\mathcal A(x,y)$ in this orthonormal basis gives
the asserted formula.
Finally, if $\omega^\alpha(x)=0$ for every $\alpha$, then
$\mathcal A(x,y)=0$ for every $y$, so $x=0$. The linear map
$$
x\longmapsto\bigl(\omega^1(x),\ldots,\omega^n(x)\bigr)
$$
is therefore a bijection from $\R^n$ onto $\R^n$.
Consequently, the $\omega^\alpha$ form a basis of the dual of $\E^n$.
\subsection{A classification of flat forms in dimensions $n\leq3$}
When the domain has dimension $n\leq3$, Cartan's second result, stated in the introduction, classifies all flat forms on $\R^n$, regardless of the dimension of the span of their image:
\begin{theorem} \label{T2}
Let ${\cal A}:\R^n\times\R^n\to\E^p$, $n\leq3$, be a flat symmetric bilinear form. Then there exist an orthogonal $p\times p$ matrix $(a^\alpha_\beta)$ and $p$ linear forms $\omega^1,\ldots,\omega^p$ such that, for every $\alpha\in\{1,\ldots,p\}$,
$$A^{\alpha} = \sum_{\beta=1}^{p}a_{\beta}^{{\alpha}}\omega^{\beta} \otimes \omega^{\beta}.$$
\end{theorem}
\begin{corollary}\label{201}
Let ${\cal A}:\R^n\times\R^n\to\E^p$, $1 \leq n\leq3$, be a flat symmetric bilinear form  of rank $r$. 
 Then there exist an orthonormal basis
$\nu^1,\ldots,\nu^p$ of $\E^p$ and $r$ linear forms $\omega^1,\ldots,\omega^r$ of rank $1$ such that $${\cal A}=\sum_{\alpha=1}^r\omega^\alpha\otimes\omega^\alpha\otimes\nu^\alpha.$$
\end{corollary}
\vspace{10pt}
Let us examine the possible values of $n$ and the rank $r$. The case $r=1$ is trivial by observation  \ref{Par212} of section{Par21}. 
\begin{enumerate}
\item
Suppose that $n=2$. Since $\dim\Sym(\R^2)=3$, we have $r\leq3$. If $r\leq2$, Theorem~\ref{T1} gives the conclusion. If $r=3$, then $F=\Sym(\R^2)$, so every rank-one form belongs to $F$. After changing the basis, we may assume that $A_1$ has rank $1$. By observation~\ref{Par213}, ${\cal A}'=\{A_2\}$ is flat. Theorem~\ref{T1}, with $n=2$, completes the argument.
\item
Suppose that $n=3$. Since $\dim\Sym(\R^3)=6$, we must consider $r\in\{1,2,3,4,5,6\}$.
	\begin{enumerate}
		\item
If $r\leq3$, Theorem~\ref{T1} gives the conclusion.
		\item
Consider the case $r=4$.
If $\ker\mathcal A\neq\{0\}$, passing to the quotient $\R^3/\ker\mathcal A$ reduces the problem to dimension at most $2$, which has already been treated. If $\ker\mathcal A=\{0\}$, Proposition~\ref{P2} asserts that $F$ contains a positive definite form, which can be represented by the identity matrix $I$ in a suitable basis.
We use the following lemma (see, for example, \cite{Brickman1961}):
\vspace{10pt}
\begin{lemma} \label{L2}
Let $q_1,q_2$ be nonzero quadratic forms of trace zero on a vector space of dimension $n\geq3$. Then $q_1$ and $q_2$ have a common isotropic vector: there exists a nonzero vector $v$ such that $q_1(v)=q_2(v)=0$.
\end{lemma}
\vspace{10pt}
\noindent Apply this lemma to the two-dimensional space $F^\perp$.
Choose a basis $B^1,B^2$ of $F^\perp$. Then $\langle I,B^1\rangle=\tr B^1=0$ and $\langle I,B^2\rangle=\tr B^2=0$. Since $n>2$, Lemma~\ref{L2} implies that $B^1,B^2$ have a common isotropic vector $v$. Thus $B^1(v,v)=B^2(v,v)=0$, or equivalently $\langle B^1,v\otimes v\rangle=\langle B^2,v\otimes v\rangle=0$. Identifying $\R^3$ with its dual, for instance by its standard inner product, gives $v\otimes v\in F$. Hence $F$ contains the rank-one form $v\otimes v$. Choose a basis $(C^1,C^2,C^3,C^4)$ of $F$ with $C^4=v\otimes v$. By Section~\ref{Par21}, observation~\ref{Par213}, the family $(C^1,C^2,C^3)$ is flat. This reduces the problem to the preceding case, where Theorem~\ref{T1} applies.
		\item
If $r=5$, then $\dim F^\perp=1$. The same proof applies to a generator of $F^\perp$ admitting an isotropic vector.
		\item 
If $r=6$, then $F=\Sym(\R^3)$, and the conclusion is trivial since $F$ contains all rank-one forms.
\end{enumerate}
\end{enumerate}
Thus, in every case, $F$ admits a basis of rank-one forms. This proves Theorem~\ref{T2}. Its corollary is immediate: as in the proof of Corollary~\ref{C1}, it remains only to check that the change from the standard basis to the basis $\nu^\alpha$ is orthogonal.
\noindent
\begin{remark}
The key argument underlying Theorem~\ref{T2} is the existence of a common isotropic vector for two quadratic forms on $\R^n$ when $n>3$.
\begin{enumerate}
\item
This result fails for two quadratic forms on $\R^2$. Indeed, the only solution of $z^2=0$ in $\C$ is $z=0$. Writing $z=x+\mathrm{i}y$, this means that the only solution of
$x_1^2-x_2^2=0$, $x_1x_2=0$ is $(0,0)$.
\item
One might hope for an analogue of Theorem~\ref{T2} for three trace-zero quadratic forms. Unfortunately, this is false, suggesting that different arguments are needed in higher dimensions.
For example, in $\R^3$, the three linearly independent trace-zero forms $q_1=x_1^2-x_2^2$, $q_2=x_2^2-x_3^2$, and $q_3=x_1x_2+x_2x_3+x_3x_1$ have no common isotropic vector.
The same holds in $\R^4$ for the three quadratic forms
$
q_1(x) =x_1^2 + x_2^2 - x_3^2 - x_4^2,
q_2(x) = 2(x_1x_3 + x_2x_4),
q_3(x) = 2(x_1x_4 - x_2x_3).
$
\end{enumerate}
\end{remark}
Cartan stopped  here this classification.
The first case not covered in \cite{Cartan1919} is that of a flat form of rank $r\geq5$ on $\R^4$. This is the subject of Section~\ref{Par5}.

\section{The endomorphism ${\cal C}_{\cal A}$ of $\Sym(\R^n)$}\label{Par22}
Let ${\cal A}=(A^1,\ldots,A^p)$ be a flat on $\R^n$ form with values in $\E^p$. Our main tool in the following sections is the endomorphism ${\cal C}_{\cal A}:\Sym(\R^n)\to\Sym(\R^n)$ defined by
\begin{equation}\label{E117}
{\cal C}_{\cal A}(M)=\sum_{{\alpha}=1}^{p}
\langle A^{\alpha},M \rangle A^{\alpha}
=\sum_{{\alpha}=1}^{p}\tr(A^{\alpha} M) A^{\alpha} ,
\end{equation}
for every $M\in\Sym(\R^n)$.
Using the totally symmetric tensor $T_{\cal A}$ defined in observation~\ref{Par214} of Section~\ref{Par21}, we obtain
$$
\bigl(\mathcal C_{\cal A}(M)\bigr)_{ij}
=\sum_{k,l=1}^n T_{{\cal A}_{ijkl}}M_{kl}.
$$
(The last two indices of $T$ are contracted with those
of $M$.)
Clearly, ${\cal C}_{\cal A}$ is self-adjoint,
since
$\langle{\cal C}_{\cal A}(M),N\rangle
=\sum_{{\alpha}=1}^{n}
\langle A^{\alpha},M \rangle \langle A^{\alpha},N \rangle$
and positive semidefinite, since
\begin{equation} \label{E12}
\langle{\cal C}_{\cal A}(M),M\rangle
=\sum_{{\alpha}=1}^{p}\langle A^{\alpha},M\rangle^2\geq 0.
\end{equation}
Hence
$\ker\mathcal C_{\cal A}=F^\perp$ and
$\rank\mathcal C_{\cal A}=\dim F=r$,
which implies
$\dim \ker {\cal C}_{\cal A}=\frac{n(n+1)}2-r$.
We summarize this as follows.
\begin{proposition} \label{P3}
The endomorphism ${\cal C}_{\cal A}$ of $\Sym(\R^n)$ is self-adjoint and positive semidefinite. Its rank is $r$, and its kernel has dimension $\frac{n(n+1)}2-r$.
\end{proposition}
As explained in Appendix~\ref{Annexe}, we may consider the functional $\ell_{T_{\cal A}}$, abbreviated to $\ell_{\cal A}$, defined on $\mathcal P_4$ by linear extension of
\begin{equation} \label{E957}
\ell_{\cal A}(x_ix_jx_kx_l) = T_{{\cal A}ijkl}.
\end{equation}
Identifying a symmetric matrix $M$ with its associated quadratic form $q_M$, and hence with a homogeneous polynomial of degree $2$, we may associate with $\ell_{\cal A}$ a symmetric bilinear form $B_{\ell_{\cal A}}$ on the space of quadratic forms by setting, for all $M,N\in\Sym(\R^n)$,
\begin{equation}\label{E14}
B_{\ell_{\cal A}}(q_Mq_N)=\ell_{\cal A}(q_Mq_N)
=\langle\mathcal C_{\cal A}(M),N\rangle
=\sum_{{\alpha}=1}^{p}
\langle A^{\alpha},M\rangle\langle A^{\alpha},N\rangle.
\end{equation}
The preceding discussion yields the following lemma.
\begin{lemma} \label{L1}
For all $M,N\in\Sym(\R^n)$, we have
\begin{enumerate}
\item
$\ell_{\cal A}(q_M^2) =\sum_{{\alpha}=1}^{p}\langle A^{\alpha},M\rangle^2\geq0$,
\item
$\ker (B_{\ell_{\cal A}})=\{q_M : M  \in \ker {\cal C}_{\cal A}\}
=\{q_M : M\in F^\perp\}$,
\item
$\rank B_{\ell_{\cal A}} = \dim F$,
\end{enumerate}
where $q_M$ and $q_N$ are the quadratic forms associated with $M$ and $N$, respectively.
\end{lemma}
We thus obtain the following corollary.
\begin{corollary} \label{C21}
The linear functional $\ell_{\cal A}$ is nonnegative on the cone of squares of
quadratic forms.
\end{corollary}
\section{Flat forms of rank $r=5$ on $\R^4$}\label{Par5}
Since $\dim\Sym(\R^4)=10$, we shall consider flat forms of ranks $5,6,7,9,10$. We leave aside the case $r=8$, for which the existence of a basis of rank-one forms spanning $\cal A$ remains conjectural.
\subsection{Flat forms of rank $r\leq5$ on $\R^4$}
We shall prove the following general result.
\begin{theorem} \label{T6}
Let ${\cal A}:\R^4\times\R^4\to\E^p$ be a flat symmetric bilinear form such that $\dim \Vect(\Im {\cal A})=r$.
If $1\leq r\leq5$, then $\Vect(\Im {\cal A})$ admits a basis of $r$ rank-one linear forms.
\end{theorem}
In other words, there exist nonzero linear forms $\omega^1,\ldots,\omega^r$
such that
 $\Vect(\omega^1 \otimes \omega^1, \ldots , \omega^r \otimes \omega^r) =\Vect (\Im {\cal A}).$
\begin{corollary} \label{C6}
Let ${\cal A}:\R^4\times\R^4\to\E^p$ be a flat symmetric bilinear form such that $\dim\Vect(\Im\mathcal A)=r$.
If $1\leq r\leq5$, there exist an \emph{orthonormal} basis $(\nu^1,\ldots,\nu^p)$ of $\E^p$ and $r$ rank-one linear forms $\omega^1,\ldots,\omega^r$ such that
$${\cal A} = \sum_{{\alpha} =1}^{r} \omega^{\alpha} \otimes \omega^{\alpha} \otimes \nu^{\alpha} .$$  
\end{corollary}
\vspace{10pt}
\noindent\textbf{Proof of Theorem~\ref{T6}.}
It suffices to prove the theorem for $r=5$; for $r<5$, the results of the preceding sections apply.
Since $\cal A$ is flat, $T_{\cal A}$ is totally symmetric and ${\cal C}_{\cal A}$ belongs to the cone $\mathcal K$ of positive semidefinite self-adjoint endomorphisms of $\Sym(\R^4)$ arising from totally symmetric tensors, as in~\ref{E18}. Using the decomposition~\ref{E20}, we obtain
$${\cal C}_{\cal A} =  {\cal C}^1 + \ldots +{\cal C}^{\alpha}+ \ldots + {\cal C}^k,$$
where each ${\cal C}^\alpha$ is extremal in $\mathcal K$. We may assume that every ${\cal C}^\alpha$ is nonzero. By Proposition~\ref{P3}, ${\cal C}_{\cal A}$ has rank $r=5$, so each ${\cal C}^\alpha$ has rank at most $5$.
We now use Proposition~\ref{P55}, which states that, in dimension $n=4$, the bilinear form associated with an extreme ray of the dual cone of sums of squares of quadratic forms has rank $1$ or $6$. Under the hypotheses of Theorem~\ref{T6}, each ${\cal C}^\alpha$ generates an extreme ray of rank at most $5$. Hence each ${\cal C}^\alpha$ must be a symmetric endomorphism of $\Sym(\R^4)$ of rank at most $1$. There is therefore a rank-one matrix $S^\alpha\in\Sym(\R^4)$ such that
\begin{equation}\label{E101} 
{\cal C}^{\alpha}(M)= \langle S^{\alpha}, M \rangle S^{\alpha},
\end{equation}
for every $M\in\Sym(\R^n)$.
We also know that, for every $\alpha$,
\begin{equation} \label{E102} 
\Im {\cal C}^{\alpha} \subset \Im {\cal C}_{\cal A},
\end{equation} 
so $\Vect(\Im\mathcal A)=\Vect(S^1,\ldots,S^k)$. Extracting a basis of five elements from $S^1,\ldots,S^k$ completes the proof.
\begin{remark}\label{91}
Once ${\cal C}_{\cal A}$ has been decomposed into a sum of rank-one endomorphisms, one could alternatively use observation~\ref{Par213} of Section~\ref{Par21} to reduce to $r\leq4$ and apply Theorems~\ref{T1} and~\ref{T2} directly.
\end{remark}
\noindent\textbf{Proof of Corollary~\ref{C6}.} It remains to find an \emph{orthonormal} basis $(\nu_1,\ldots,\nu_r)$ of the codomain of $\cal A$ such that
$${\cal A}  =  \sum_{{\alpha}=1}^{r}\omega^{\alpha} \otimes \omega^{\alpha} \otimes \nu_{\alpha}.$$
If $r\leq4$, this follows from Theorems~\ref{T1} and~\ref{T2}. If $r=5$, Theorem~\ref{T6} provides a unit vector $\nu_5$ such that $\langle\mathcal A(\cdot,\cdot),\nu_5\rangle$ has rank $1$. This component can be ``eliminated'' by observation~\ref{Par213} of Section \ref{Par21}, reducing to rank $4$. Theorems~\ref{T1} and~\ref{T2} then give the conclusion.
\subsection{Flat forms of rank $6$ or $7$ on $\R^4$}
The situation changes completely when $5$ is replaced by $6$ or $7$.
\subsubsection{Flat forms on $\R^4$ of rank $6$}
In his doctoral thesis \cite{Zafindratafa1986} (see also \cite{MorvanZafindratafa1987}), G.~Zafindratafa showed that Theorems~\ref{T6} and~\ref{T7} do not extend to the case where the domain has dimension $4$ and $\dim\Vect(\Im\mathcal A)=6$. He exhibited the following six $4\times4$ matrices. They are linearly independent, their kernels have trivial intersection, they form a flat family, and no linear combination of them has rank $1$.
Here is the example:
{\tiny
$$
A^1=\begin{pmatrix}
1 &0 &0 &0 \\
 0 &1& 0 &0 \\
0 &0& 1& 1 \\
0& 0 &1 &1 \\
\end{pmatrix},
A^2=\begin{pmatrix}
0&1&0&0 \\
1&0&0&0 \\
0&0&0&0 \\
0&0&0&0 \\
\end{pmatrix},
A^3=\begin{pmatrix}
0&0&0&0 \\
0&0&1&0 \\
0&1&0&0 \\
0&0&0&0 \\
\end{pmatrix},
$$
$$
A^4=\begin{pmatrix}
0&0&0&0 \\
0&1&0&0 \\
0&0&0&-1\\
0&0&-1&-1 \\
\end{pmatrix},
A^5=\begin{pmatrix}
0&0&1&1\\
0&0&0&0\\
1&0&0&0\\
1&0&0&0\\
\end{pmatrix}
A^6=\begin{pmatrix}
0&0&0&0\\
0&0&0&0\\
0&0&1&0\\
0&0&0&1\\
\end{pmatrix}.
$$
}
Thus no general result of this kind can be expected for $n=4$.
\subsubsection{Flat forms on $\R^4$ of rank $7$}
The conclusion is the same when $r=7$: in general, $F$ need not admit a basis of rank-one forms, as the following example shows. Consider the seven matrices
{\tiny
$$A^1=\begin{pmatrix}
\frac{4}{3}&0&0&0\\
0&\frac{1}{3}&0&0\\
0&0&\frac{1}{3}&0\\
0&0&0&\frac{1}{2}\\
\end{pmatrix},
A^2=\begin{pmatrix}
\frac{1}{3}&0&0&0\\
0&\frac{4}{3}&0&0\\
0&0&\frac{1}{3}&0\\
0&0&0&\frac{1}{2}\\
\end{pmatrix},
A^3=\begin{pmatrix}
\frac{1}{3}&0&0&0\\
0&\frac{1}{3}&0&0\\
0&0&\frac{4}{3}&0\\
0&0&0&\frac{1}{2}\\
\end{pmatrix},
$$
$$
A^4=\begin{pmatrix}
0&1&0&0\\
1&0&0&0\\
0&0&0&1\\
0&0&1&0\\
\end{pmatrix},
A^5=\begin{pmatrix}
0&0&1&0\\
0&0&0&1\\
1&0&0&0\\
0&1&0&0\\
\end{pmatrix},
A^6=\begin{pmatrix}
0&0&0&1\\
0&0&1&0\\
0&1&0&0\\
1&0&0&0\\
\end{pmatrix},
A^7=\begin{pmatrix}
1&1&0&0\\
1&1&0&0\\
0&0&0&0\\
0&0&0&0\\
\end{pmatrix}.
$$
}
A lengthy calculation shows that these matrices are linearly independent and form a flat family. Nevertheless, the following weaker general result holds; it can also be deduced from \cite{Blekherman2015}.
\begin{proposition}\label{P98}
Let ${\cal A}=(A^1, \ldots, A^p):\R^4\times\R^4\to\E^p$ be a flat symmetric bilinear form of rank $7$.
Then $F$ contains at least one rank-one form.
\end{proposition}
\noindent\textbf{Proof of Proposition~\ref{P98}.} We again use a result of Blekherman (corollary~1.4 of \cite{Blekherman2012}), stated in our setting as follows: \emph{a quartic form strictly positive on $\R^{4*}$ located on the boundary of  $\Sigma_{4,4}$ cannot be expressed as a sum of three squares of quadratic forms}. The orthogonal complement of $F$ has dimension $3$. Choose three linearly independent symmetric bilinear forms in $F^\perp$, with associated quadratic forms $q_1,q_2,q_3$. The quartic functional $\ell_{\cal A}$ associated with $\cal A$ is positive, and its kernel is $F^\perp$, as in Lemma~\ref{L1}. The quartic form $q=q_1^2+q_2^2+q_3^2$ therefore lies on the boundary of $\Sigma_{4,4}$, the space of symmetric bilinear forms on $\R^4$ that are sums of squares of quadratic forms. If $q_1,q_2,q_3$ had no common isotropic vector, $q$ would be positive definite, contradicting \cite{Blekherman2012}. Thus there is a common isotropic vector $x_0$. Identifying $\R^4$ with its dual through an inner product $\langle\cdot,\cdot\rangle$, the linear form $\omega_0=\langle x_0,\cdot\rangle$ satisfies $\omega_0\otimes\omega_0\in F$, and this form has rank $1$.

\subsection{Flat forms on $\R^4$ of rank $8$, $9$, or $10$}
\subsubsection{The case $r=9$}
We obtain the following theorem. Its corollary asserts the existence of a basis of $\R^9$ satisfying~\ref{E96}, but in general this basis cannot be chosen orthonormal.
\begin{theorem} \label{T9}
Let $(A^1,\ldots,A^9)$ be a linearly independent flat family of real symmetric bilinear forms on $\R^4$.
Then $F$ admits a basis of positive semidefinite rank-one symmetric bilinear forms.
\end{theorem}
In other words, there exist nonzero linear forms
$\omega^1,\ldots,\omega^9\in\R^4$ such that
$$
 F=\Vect(\omega^1\otimes \omega^1,\ldots, \omega^p \otimes \omega^p).
$$
\begin{corollary} \label{C69}
Let ${\cal A}=(A^1,\ldots,A^p):\R^4\times\R^4\to\E^p$ be a flat symmetric bilinear form.
Set $F=\Vect(A^1,\ldots,A^p)$ and assume that $\dim F=9$.
Then there exist a basis $(\nu^1,\ldots,\nu^9)$ of $\E^9$ and nine positive semidefinite rank-one forms $\omega^1,\ldots,\omega^9$ such that
\begin{equation}\label{E96}
{\cal A} = \sum_{{\alpha} =1}^{9} \omega^{\alpha} \otimes \omega^{\alpha} \otimes \nu^{\alpha} .
\end{equation} 
\end{corollary}

\noindent\textbf{Proof of Theorem~\ref{T9}~-} Here $F$ is a hyperplane in $\Sym(\R^4)$. Let $J$ be a bilinear symmetric matrix spanning $F^\perp$. We can assume, with a suitable  frame of $\R^4$, that $H$ is a diagonal matrix $H=\diag ( h_1,h_2,h_3,h_4)$. Using the fact that $\cal A$ is flat, have
$$\sum_{\alpha=1}^9 A^\alpha_{ii}A^\alpha_{jj}
=\sum_{\alpha=1}^9 (A^\alpha_{ij})^2.$$
Then for every $j \in \{1, \ldots , 4\}$, one gets, using the fact that $H$ is orthogonal to each $A^\alpha$,
\begin{equation}\label{E7221}
\sum_{i=1}^4 h_i\sum_{\alpha=1}^9(A^\alpha_{ij})^2
=\sum_{\alpha=1}^9 A^\alpha_{jj}\,\langle H,A^\alpha\rangle
=0.
\end{equation}
\begin{lemma}\label{L869}
The form $H$ is indefinite.
\end{lemma}
Indeed, if $H$ was semidefinite (for instance positive), all $h_i$ should be greater or equal to $0$. Since $H$ is not null, there should exist $k$ such that $h_k >0$. From equality {E7221}, we deduce that for every $\alpha$ and every  $j$, $A_{kj}^\alpha = 0$. This would implies that $F$ is include in a vector subspace of dimension $6$ in $\Sym(\R^n)$. This is exclude since $\dim F=9$.
Now, we can apply  the following lemma \cite{Elton2010}~:
\begin{lemma} \label{Isot}
 Let  $q,r$ be two quadratic forms on $\R^n$. Suppose that   there exist $u,v \in \R^n$ such that $q(u)>0$ and $q(v)<0$. If
$(q(x)=0\ \Longrightarrow\ r(x)=0)$,
then there exists a real number $a$ such that
$\forall x \in \R^n,  r(x)=aq(x)$.
\end{lemma}

Let $G=\Vect\{x\otimes x: H(x,x)=0\}$. Clearly $H \in G^\perp$. Conversly, let $B\in G^\perp$. For every $x$ such that  $H(x,x)=0$, one has $B(x,x)=0$. Since $H$ is indefinite, Lemma \ref{Isot} implies that  $q_B=a q_H$  for some real number $a$. We conclude that $G=F$.

\vspace{10pt}

\begin{remark}
Corollary~\ref{C69} does not imply the existence of an \emph{orthonormal} basis of $\Vect(\Im\mathcal A)$ satisfying~\ref{E96}. For an example, consider the following nine matrices:
{\tiny $$A^1=\begin{pmatrix}
\frac{4}{3}&0&0&0\\
0&\frac{1}{3}&0&0\\
0&0&\frac{1}{3}&0\\
0&0&0&\frac{1}{2}\\
\end{pmatrix},
A^2=\begin{pmatrix}
\frac{1}{3}&0&0&0\\
0&\frac{4}{3}&0&0\\
0&0&\frac{1}{3}&0\\
0&0&0&\frac{1}{2}\\
\end{pmatrix},
A^3=\begin{pmatrix}
\frac{1}{3}&0&0&0\\
0&\frac{1}{3}&0&0\\
0&0&\frac{4}{3}&0\\
0&0&0&\frac{1}{2}\\
\end{pmatrix},
$$
$$
A^4=\begin{pmatrix}
0&1&0&0\\
1&0&0&0\\
0&0&0&1\\
0&0&1&0\\
\end{pmatrix},
A^5=\begin{pmatrix}
0&0&1&0\\
0&0&0&1\\
1&0&0&0\\
0&1&0&0\\
\end{pmatrix},
A^6=\begin{pmatrix}
0&0&0&1\\
0&0&1&0\\
0&1&0&0\\
1&0&0&0\\
\end{pmatrix},
$$
$$
A^7=\frac{1}{4}\begin{pmatrix}
1&1&0&0\\
1&1&0&0\\
0&0&0&0\\
0&0&0&0\\
\end{pmatrix},
A^8=\frac{1}{4}\begin{pmatrix}
1&0&0&0\\
0&0&0&0\\
0&0&0&0\\
0&0&0&0\\
\end{pmatrix},
A^9=\frac{1}{4}\begin{pmatrix}
1&0&1&0\\
0&0&0&0\\
1&0&1&0\\
0&0&0&0\\
\end{pmatrix}.
$$
}
\end{remark}
One checks that $A^1,\ldots,A^9$ are linearly independent, that their kernels have trivial intersection, and that ${\cal A}=(A^1,\ldots,A^9)$ is flat. However, no orthonormal basis of $\Vect(\Im\mathcal A)$ satisfies~\ref{E96}. To see this, it suffices to find a nonnegative homogeneous polynomial $f$ of degree $4$ on $\R^4$ such that $\ell_{\cal A}(f)<0$. If an orthonormal basis of $\R^9$ satisfying equality~\ref{E96} existed, then $\ell_{\cal A}(f)\geq0$, a contradiction. In our case, 
the quartic form
$$f(x)=x_1^2x_2^2   +  x_2^2x_3^2 + x_3^2x_1^2 + x_4^4  - 4x_1x_2x_3x_4. $$       
satisfies these conditions: one checks that
\begin{equation}\label{E934}
f(x)  \geq 0, \forall x \in \R^4, \ \ \ell_{\cal A}(f) =-\frac{1}{8}.
\end{equation} 
\subsubsection{The case $r=8$}
The situation is more delicate. Consider the first eight matrices ${\cal A}=\{A^1,\ldots,A^8\}$ from the preceding family (the ninth matrix has rank $1$). They still form a linearly independent flat family with trivial common kernel, and satisfy $\ell_{\cal A}(f)<-\frac{3}{16}$. Thus no orthonormal basis $\nu^\alpha$ of $\E^8$ can satisfy
\begin{equation}\label{E961}
{\cal A} = \sum_{{\alpha} =1}^{8} \omega^{\alpha} \otimes \omega^{\alpha} \otimes \nu^{\alpha} .
\end{equation} 
However, we do not have a complete proof of an analogue of Theorem~\ref{T9} without additional hypotheses. For example, assuming that $F^\perp$ is spanned by two simultaneously diagonalizable quadratic forms makes the proof very simple.

\subsubsection{The case $r=10$} The conclusion is the same as for $r=9$. 
Consider the $10$ following matrices~:
{\tiny
$$A^1=\begin{pmatrix}
\frac{4}{3}&0&0&0\\
0&\frac{1}{3}&0&0\\
0&0&\frac{1}{3}&0\\
0&0&0&\frac{1}{2}\\
\end{pmatrix},
A^2=\begin{pmatrix}
\frac{1}{3}&0&0&0\\
0&\frac{4}{3}&0&0\\
0&0&\frac{1}{3}&0\\
0&0&0&\frac{1}{2}\\
\end{pmatrix},
A^3=\begin{pmatrix}
\frac{1}{3}&0&0&0\\
0&\frac{1}{3}&0&0\\
0&0&\frac{4}{3}&0\\
0&0&0&\frac{1}{2}\\
\end{pmatrix},
$$
$$
A^4=\begin{pmatrix}
0&1&0&0\\
1&0&0&0\\
0&0&0&1\\
0&0&1&0\\
\end{pmatrix},
A^5=\begin{pmatrix}
0&0&1&0\\
0&0&0&1\\
1&0&0&0\\
0&1&0&0\\
\end{pmatrix},
A^6=\begin{pmatrix}
0&0&0&1\\
0&0&1&0\\
0&1&0&0\\
1&0&0&0\\
\end{pmatrix},
$$
$$
A^7=\frac{1}{4}\begin{pmatrix}
1&1&0&0\\
1&1&0&0\\
0&0&0&0\\
0&0&0&0\\
\end{pmatrix},
A^8=\frac{1}{4}\begin{pmatrix}
1&0&0&0\\
0&0&0&0\\
0&0&0&0\\
0&0&0&0\\
\end{pmatrix},
A^9=\frac{1}{4}\begin{pmatrix}
1&0&1&0\\
0&0&0&0\\
1&0&1&0\\
0&0&0&0\\
\end{pmatrix},
A^{10}=\frac{1}{4}\begin{pmatrix}
1&0&0&1\\
0&0&0&0\\
0&0&0&0\\
1&0&0&1\\
\end{pmatrix}.
$$
}

It is a basis of $F=\Sym(\R^4)$ consisting of ten positive semi-definite rank-one forms. These form are a flat family $\cal A$ satisfying all the required hypotheses. But $\ell_{\cal A}=-\frac{1}{16}$. Thus there exist a basis $(\nu^1,\ldots,\nu^p)$ of $\E^{10}$ and ten rank-one forms $\omega^1,\ldots,\omega^{10}$ such that
\begin{equation}\label{E961ten}
{\cal A} = \sum_{{\alpha} =1}^{10} \omega^{\alpha} \otimes \omega^{\alpha} \otimes \nu^{\alpha} .
\end{equation} 
However,
an \emph{orthonormal} basis satisfying~\ref{E961ten}
does not exist.
\vspace{10pt}
The following table summarizes this section for flat forms $\mathcal A$ on $\R^4$ with $\ker\mathcal A=\{0\}$.
\vspace{10pt}
\begin{itemize}[label={}]
\item
$\mathbf{(A)}$ means that $F$ contains a positive semidefinite rank-one symmetric bilinear form.
\item
$\mathbf{(B)}$ means that $F$ admits a basis of positive semidefinite rank-one symmetric bilinear forms.
\item
$\mathbf{(O)}$ means that there exist an \emph{orthonormal} basis
$(\epsilon^1,\ldots,\epsilon^r)$ of $\Vect(\Im\mathcal A)$ and
nonzero linear forms $\omega^1,\ldots,\omega^r$ such that
$\mathcal A=\sum_{\alpha=1}^r
\omega^\alpha\otimes\omega^\alpha\otimes\epsilon^\alpha
$.
\item
\textbf{``Not in general''} means that counterexamples exist.
\end{itemize} 
\vspace{10pt}
\begin{center}
\renewcommand{\arraystretch}{1.16}
\begin{tabular}{cccccc}
\toprule
$n$ & $r=\dim F$ & Compatibility & \textbf{(A)} & \textbf{(B)} & \textbf{(O)} \\
\midrule
$4$ & $1$ & Impossible & --- & --- & --- \\
$4$ & $2$ & Impossible &---& --- & --- \\
$4$ & $3$ & Impossible &---& --- & --- \\
$4$ & $4$ & Possible &\oui& \oui & \oui \\
$4$ & $5$ & Possible &\oui& \oui & \oui \\
$4$ & $6$ & Possible &\non& \non & \non \\
$4$ & $7$ & Possible &\oui& \non & \non \\
$4$ & $8$ & Possible &\oui& \oui & ? \\
$4$ & $9$ & Possible &\oui& \oui & \non \\
$4$ & $10$ & Possible &\oui& \oui & \non \\
$4$ & $r>10$ & Impossible& --- & --- & --- \\
\bottomrule
\end{tabular}
\end{center}
\smallskip
\section{Flat forms on $\R^n$ of rank $n+1$}\label{Par606}
This section partially generalizes Theorem~\ref{T6} to real symmetric bilinear forms $\cal A$ on $\R^n$ with values in $\R^{n+1}$. The codomain thus has dimension one greater than in Theorem~\ref{T1}. This is only a partial generalization, since we require a trivial kernel and maximal rank, conditions not imposed in Theorem~\ref{T6}.
Identifying matrices with the symmetric bilinear forms
they represent, as usual, we obtain the following theorem.
\begin{theorem}\label{T7}
  Let
$A^1,\ldots,A^{n+1} \in\Sym(\mathbb R^n)$, $n\geq4$,
be a linearly independent flat family of real symmetric matrices satisfying
$
\bigcap_{{\alpha}=1}^{n+1}\ker A^{\alpha}=\{0\}.
$
Then there exist an orthogonal matrix
$(a^{\beta}_{\alpha})\in O(n+1)$ and
nonzero linear forms $\omega^1,\ldots,\omega^{n+1}\in(\R^n)^*$
such that 

$$
\sum_{\alpha=1}^{n+1}a^\beta_\alpha A^\alpha
=\omega^\beta\otimes\omega^\beta,
\qquad \beta=1,\ldots,n+1.
$$
Consequently, $(\omega^1 \otimes\omega^1, \ldots  , \omega^{n+1} \otimes\omega^{n+1})$ is a basis of
$
F=\operatorname{span}(A^1,\ldots,A^{n+1})
$
consisting of positive semidefinite rank-one forms.
\end{theorem}

\begin{corollary} \label{C72}
Let 
$
\mathcal A:\mathbb R^n\times\mathbb R^n\longrightarrow\mathbb E^{n+1}
, \ (n\geq 4)$,
be a flat symmetric bilinear form with values in $\E^{n+1}$.
Assume that
$
\rank\mathcal A=n+1$ and $\ker\mathcal A=\{0\}$.
Then there exist an orthonormal basis
$(\epsilon_1,\ldots,\epsilon_{n+1})$ of $\E^{n+1}$
and nonzero linear forms
$\omega_1,\ldots,\omega_{n+1}\in(\R^n)^*$ such that
$$
\mathcal A=\sum_{{\alpha}=1}^{n+1}
\omega_{\alpha}\otimes\omega_{\alpha}\otimes\epsilon_{\alpha}.
$$
\end{corollary}
\noindent\textbf{Proof of Theorem~\ref{T7}.} The proof applies a result of L.~Chen, D.~Chu, L.~Qian, and Y.~Shen (Theorems~24 and~14 of \cite{ChenChuQianShen2019}), which can be stated as follows:

\vspace{10pt}
{\it A positive totally symmetric operator $\mathcal O$ on $\R^n\times\R^n$ of rank $n+1$, whose contraction over two indices has rank $n$, admits a decomposition
${\cal O} = \sum_{s=1}^{r} v_{s}^{\otimes 4}$.}
\vspace{10pt}

With the notation of Section~\ref{Par22}, the endomorphism $\mathcal C_{\cal A}$ satisfies these conditions under the hypotheses of Theorem~\ref{T7}, by proposition \ref{P2}. Indeed, its contraction is positive definite
if and only if both $\ker\mathcal A=\{0\}$ and $\rank\mathcal C_{\cal A}=n+1$ hold. These conditions are satisfied by hypothesis, and the conclusion follows.

\vspace{10pt}

\noindent\textbf{Proof of Corollary~\ref{C72}.} Theorem~\ref{T7} provides a rank-one symmetric bilinear form in $\Vect(A^1,\ldots,A^{n+1})$, say $\omega_{n+1}\otimes\omega_{n+1}$, and a unit vector $\epsilon_{n+1}$ such that $\langle\mathcal A,\epsilon_{n+1}\rangle=\omega_{n+1}\otimes\omega_{n+1}$. By observation~\ref{Par213} of Section \ref{Par21}, the form $\mathcal B=\mathcal A-\omega_{n+1}\otimes\omega_{n+1}\otimes\epsilon_{n+1}$ is still flat, has rank $n$, and takes values in $\epsilon_{n+1}^\perp$. Moreover, $\ker\mathcal B\cap\ker\omega_{n+1}=\{0\}$, so $\dim\ker\mathcal B\leq1$. We may therefore apply Theorem~\ref{T1} directly if $\ker\mathcal B=\{0\}$, or pass to the quotient by this kernel to reduce the dimensions by one. An immediate induction completes the proof.
\section{Further directions}\label{Persp}
As recalled in the introduction, this study belongs to the theory of submanifolds of Riemannian manifolds. \emph{Classifying all isometric immersions of one Riemannian manifold into another} is a long-term objective. Here we have considered a pointwise approach under zero-curvature assumptions and dimension restrictions that are nevertheless less stringent than Cartan's. This amounts to studying the intersections of the Veronese cone in $\Sym(\R^n)$ with certain linear subspaces. Under our hypotheses, the results impose geometric restrictions on focal directions and focal points along each fiber of the normal bundle of the submanifold (see \cite{MorvanZafindratafa1987} for the somewhat more general setting of conformally flat submanifolds). A theorem valid in arbitrary dimension and codimension under weaker curvature assumptions would be welcome; see the recent paper of J.~D.~Moore \cite{Moore2025}. This is, of course, only a first step: local results involving derivatives of the geometric objects are also needed. One must then solve the \emph{Codazzi--Ricci equations}. The preceding results are useful because these equations take particularly simple forms in the dimensions considered here. To obtain local or global results,   restrictions of the kind obtained above make the computations more accessible.
\appendix
\section{Appendix: positive semidefinite endomorphisms of $\Sym(\R^n)$}\label{Annexe}
\subsection{The cones $\Sigma_{n,4}^*$ and $\mathcal K$}
Let $\mathcal P_d=\R[x_1,\ldots,x_n]_d$ denote the space of
homogeneous polynomials of degree $d$ on $\R^n$.
We shall mainly use $\mathcal P_2$, of dimension $\frac{n(n+1)}2$, and $\mathcal P_4$, of dimension $N=\binom{n+3}{4}=\frac{n(n+1)(n+2)(n+3)}{24}$.
Let $\ell\in\mathcal P_4^*$. Identifying $\Sym(\R^n)$ with $\mathcal P_2$, the functional $\ell$ induces a symmetric bilinear form $B_\ell:\Sym(\R^n)\times\Sym(\R^n)\to\R$ by setting, for all $M,N\in\Sym(\R^n)$,
$$B_l(q_M,q_N)=\ell(q_Mq_N).$$
We first identify the dual $\mathcal P_4^*$ of $\mathcal P_4$ with a linear subspace of $\End(\Sym(\R^n))$ by constructing the linear map
$$\psi : {\cal P}_4^* \to \End ( \Sym(\R^n))$$
as follows.
Let $\ell\in\mathcal P_4^*$ be the linear functional whose tensor coefficients $\ell_{ijkl}$ in the standard basis of $\mathcal P_4$ are
$\ell_{ijkl}=x_ix_jx_kx_l$. For every $M\in\Sym(\R^n)$, we have
 $$\psi(\ell)(M)_{kl}= \sum_{ij}l_{ijkl}M_{ij}.$$
Write $\mathcal C_\ell=\psi(\ell)$.
We then have the following lemma.
\begin{lemma}
For all $M,N\in\Sym(\R^n)$,
\begin{itemize}
\item
$\langle {\cal C}_\ell (M),N \rangle = \langle {\cal C}_\ell (N),M \rangle =\ell (q_Mq_N) = B_\ell(q_M,q_N)$.
\item
$\rank B_\ell = \rank {\cal C}_\ell, \ker B_l =\{q_M : M \in \ker {\cal C}_\ell \}$.
\end{itemize}
In particular, $\mathcal C_\ell$ is self-adjoint.
\end{lemma}
The map $\psi$ is injective, but not bijective when $n>1$, because the tensor $\ell_{ijkl}$ is totally symmetric. This leads us to introduce the following linear subspace $\mathcal H$ of $\End(\Sym(\R^n))$:
$$\mathcal H=\{\mathcal C\in\End(\Sym(\R^n)):\mathcal C\text{ is totally symmetric}\}.$$
Then $\psi$ maps onto $\mathcal H$ and induces a vector-space isomorphism
$${\cal P}^*_4 \simeq {\cal H}.$$
\subsubsection{A first cone}
Identifying quadratic forms on $\R^n$ with elements of $\mathcal P_2$ as usual, introduce the convex cone $\Sigma_{n,4}\subset\mathcal P_4$ of sums of squares of quadratic forms:
\begin{equation} \label{E15}
 \Sigma_{n,4}
 =\left\{\sum_{a=1}^{s}q_a^2:s\geq1,\ q_a\in {\cal P}_2\right\}
 \subset {\cal P}_4.
\end{equation}
Its dual is the closed convex cone
\begin{equation} \label{E16}
 \Sigma_{n,4}^*
 =\{\ell\in {\cal P}_4^*:\ell(f)\geq 0, \forall f\in\Sigma_{n,4}\} \subset {\cal P}_4^*.
\end{equation}
The following lemma is immediate.
\begin{lemma} \label{L921}
$\ell\in\Sigma_{n,4}^*
 \iff \ell(q^2)\geq 0,  \forall q\in {\cal P}_2$;
\end{lemma}
\subsubsection{A second cone}
We also introduce the cone $\mathcal K$ of positive semidefinite endomorphisms $\mathcal C$ of $\Sym(\R^n)$ whose components $\mathcal C_{ijkl}$ are totally symmetric:
\begin{equation}\label{E17}
\displaystyle
 {\cal K}= \{ {\cal C} \in {\cal H}:
 \langle{\cal C}(M),M\rangle \ \geq 0,  \forall M\in \Sym(\R^n) \}.
\end{equation}
The restriction of the linear map $\psi$ to $\Sigma_{n,4}^*$ has image $\mathcal K$. We have
\begin{equation}\label{E18}
\forall \ell\in\Sigma_{n,4}^*,  \ell(q_M^2)=\langle{\cal C}_\ell(M),M\rangle \geq 0, \forall M \in \Sym(\R^n).
\end{equation}
and obtain the equivalence
$$ \ell\in \Sigma_{n,4}^* \iff {\cal C}_\ell \in {\cal K}.$$
Since $\psi$ is injective, we obtain:
\begin{proposition}\label{P67}
The linear map $\psi$ identifies the cone $\Sigma_{n,4}^*\subset\mathcal P_4^*$ with the cone $\mathcal K\subset\mathcal H$.
\end{proposition}

\subsection{Decomposition of elements of $\mathcal K$ into extreme rays}
We use the following standard property of finite-dimensional cones:
\begin{proposition}\label{P88}
Every finite-dimensional closed convex cone with a compact base is generated by its extreme rays. Thus each element of the cone is a finite linear combination, with positive coefficients, of generators of extreme rays.
\end{proposition}
Here $\mathcal K$ is closed in a finite-dimensional space, and its elements are nonnegative self-adjoint endomorphisms, identified with positive semidefinite symmetric bilinear forms. Its section $\{H\in\mathcal K:\tr H=1\}$ by an affine hyperplane is compact, since it is closed and bounded.
Consequently, every element $\mathcal C$ of this cone is a finite sum of elements lying on its extreme rays:
\begin{equation}\label{E20}
{\cal C} =  {\cal C}^1 + \ldots + {\cal C}^{\alpha}+ \ldots + {\cal C}^k,
\end{equation}
where each $\mathcal C^\alpha$ generates an extreme ray.
Moreover, since $\mathcal C$ and every $\mathcal C^\alpha$ are positive semidefinite,
\begin{enumerate}
\item
 $\bigcap_{\alpha} \ker {\cal C}^{\alpha} =\ker {\cal C}$, 
\item
$\rank\mathcal C^\alpha\leq\rank\mathcal C$ for every $\alpha$, as endomorphisms of $\Sym(\R^n)$.
\end{enumerate}

\subsection{Quartics in four variables}\label{Par243}
Restrict now to dimension $n=4$, with $\mathcal K\subset\psi(\Sigma_{4,4})^*$.
The following proposition is a consequence of the results in \cite{Blekherman2012}.
\begin{proposition}\label{P55}
Let $\mathcal C\in\mathcal K\subset\Sigma_4^*$ generate an extreme ray, and let $\ell=\psi^{-1}(\mathcal C)$ be its associated linear functional. There are only two possibilities:
\begin{enumerate}
\item
The form $B_\ell$ has rank $1$, and $\ell$ is a real evaluation: there exist $v\in\R^4\setminus\{0\}$ and $\lambda>0$ such that $\ell(f)=\lambda f(v)$ for every $f\in\mathcal P_4$.
\item
The form $B_\ell$ has rank $6$, and its kernel has dimension $4$. Moreover, this kernel is spanned by four quadratic forms with no common isotropic vector in $\C^4$.
\end{enumerate}
\end{proposition}


\begin{thebibliography}{99}
\bibitem{Cartan1919}
\textsc{Élie Cartan}, Sur les variétés de courbure constante d'un espace euclidien ou non-euclidien, \emph{Bulletin de la Société Mathématique de France}, \textsc{47} (1919), pp.~125--160.
\bibitem{Cartan1920}
\textsc{Élie Cartan}, Sur les variétés de courbure constante d'un espace euclidien ou non-euclidien, \emph{Bulletin de la Société Mathématique de France}, \textsc{48} (1920), pp.~132--208.
\bibitem{Moore1972}
\textsc{John Douglas Moore}, Isometric immersions of space forms in space forms, \emph{Pacific Journal of Mathematics}, \textsc{40} (1972), pp.~157--166.
\bibitem{Moore1977Positive}
\textsc{John Douglas Moore}, Submanifolds of constant positive curvature I, \emph{Duke Mathematical Journal}, \textsc{44} (1977), pp.~449--484.
\bibitem{Moore2025}
 \textsc{John Douglas Moore}, Submanifolds of constant negative curvature~: A generalization of Hilbert's theorem, \emph{arXiv:2512.00973v1} [math.DG] 30 Nov 2025.

\bibitem{DajczerTojeiro2019}
\textsc{Marcos Dajczer and Ruy Tojeiro}, \emph{Submanifold Theory: Beyond an Introduction}, Universitext, Springer, 2019.

\bibitem{Brickman1961}
 \textsc{Louis Brickman}, On the Field of Values of a Matrix, \emph{Proceedings of the American Mathematical Society}, \textbf{12} (1961), pp.~ 61--66.

\bibitem{ChernKuiper1952}
\textsc{S.-S. Chern and N.H. Kuiper}, Isometric imbedding of Riemannian manifolds with nonpositive curvature, \emph{Annals of Mathematics}, \textbf{56} (1952), pp.~422--430. 


\bibitem{MooreMorvan1978}
\textsc{John Douglas Moore and Jean-Marie Morvan}, Sous-variétés conformément plates de codimension quatre, \emph{Comptes rendus de l'Académie des sciences de Paris, Série A}, \textbf{287} (1978), pp.~655--657.
\bibitem{MooreMorvan2001}
\textsc{John Douglas Moore and Jean-Marie Morvan}, On isometric Lagrangian immersions, \emph{Illinois Journal of Mathematics}, \textbf{45} (2001), pp.~833--849.

\bibitem{Chen1973}
\textsc{Bang-Yen Chen}, Geometry of Submanifolds, \emph{Marcel Dekker}.


\bibitem{MorvanZafindratafa1987}
\textsc{Jean-Marie Morvan and Georges Zafindratafa}, Conformally flat submanifolds, \emph{Annales de la Faculté des sciences de Toulouse : Mathématiques}, Series~5, 
\textbf{8}, no.~3 (1986--1987), pp.~331--348.

\bibitem{Zafindratafa1986}
\textsc{Georges Zafindratafa}, Sous-variétés conformément plates d'un espace euclidien, Ph.D. thesis, Université de Provence, 20 March 1986.

\bibitem{Blekherman2012}
\textsc{Grigoriy Blekherman}, Nonnegative Polynomials and Sums of Squares, \emph{Journal of the American Mathematical Society}, \textsc{25}, no.~3 (2012), pp.~617--635.
\bibitem{Blekherman2015}
\textsc{Grigoriy Blekherman}, Positive Gorenstein Ideals, \emph{Proceedings of the American Mathematical Society}, \textsc{143} (2015), pp.~69--86.

\bibitem{ChenChuQianShen2019}
\textsc{Lin Chen, Delin Chu, Lilong Qian and Yi Shen}, Separability of completely symmetric states in a multipartite system, \emph{Physical Review A}, \textsc{99} (2019), article~032312.

\bibitem{Elton2010}
\textsc{John H. Elton}, Indefinite Quadratic Forms and the Invariance of the Interval in Special Relativity, \emph{The American Mathematical Monthly}, \textsc{117} (2010), no.~ 6, p. 540–-547.



\end{thebibliography}
\end{document}